\documentstyle[amssymb]{amsart}
\title[QKZ equation]
{A solution of the quantum Knizhnik 
Zamolodchikov equation of type $C_n$}
\author{\sc Katsuhisa Mimachi}

\date{}
\numberwithin{equation}{section}
\newtheorem{thm}{Theorem}[section]
\newtheorem{cor}[thm]{Corollary}
\newtheorem{prop}[thm]{Proposition}
\newtheorem{df}[thm]{Definition}
\newtheorem{lem}[thm]{Lemma}

\begin{document}
\maketitle
\begin{abstract}
We construct a solution of Cherednik's quantum Knizhnik
Zamolodchikov equation associated with the root system
of type $C_n\,.$  This solution is given in terms of a 
restriction of a $q$-Jordan-Pochhammer integral. 
As its applicaton, we give an explicit expression
of a special case of the Macdonald polynomial of 
the $C_n$ type.
Finally we explain the connection with the 
representation of the Hecke algebra. 
\end{abstract}
\section{Introduction}
We study the quantum Knizhnik Zamolodchikov (QKZ) 
equation (\cite{C1}) associated with the root 
system of type $C_n.$ 
A solution to this equation is found by means 
of a restriction of the $q$-Jordan-Pochhammer 
integral. 

A solution of the QKZ equation of type $A_{n-1}$ 
is given in \cite{Mi1}. 
Since the appearance of that work, however, 
there has been no progress in the study of 
the QKZ equation for other types of root systems
with regard to the determination of solutions. 
This paper is devoted to such a task.

To construct our solution, we exploit a 
family of rational functions which
would correspond to a basis of the $q$ 
de Rham cohomology attached to the 
integrand.  This turns out to be a natural 
basis for the representation of the 
Hecke algebra $H(W)$ through the 
Lusztig operator $T_i\,.$

Next, as a byproduct of our investigation, 
we obtain an integral representaion of the
special case of an eigenfunction 
associated with the Macdonald operator of 
the $C_n$ type.
In particular, it is seen that, taking a 
suitable cycle, a restriction of the 
$q$-Jordan-Pochhammer integral
expresses the Macdonald polynomial of the
$C_n$ type parametrized by
the partition $(\lambda,0,\ldots,0).$ 
This integral leads to a more explicit
expression.\smallskip 

We believe that the present paper represents
a first step toward understanding the $BC_n$ 
type QKZ equation and the $BC_n$ type 
Macdonald polynomial. 
It is noteworthy that even in 
the classical ($q$=1) case was not previously 
known that such an integral gives spherical 
functions associated with the root system $C_n.$
For related works on $BC_n$ type spherical 
functions, we refer the reader to \cite{DG} 
and references therein.
\smallskip

Throughout this paper, $q$ is regarded as a real 
number satisfying $0\le q <1\,.$ 

\section{QKZ equation of type $C_n$}
We first give a review of the QKZ equation 
associated with the root system of type $C_n$ 
for the reader's convenience, following 
Cherednik \cite{C1} and Kato \cite{Kato}.\smallskip

Let $E=\oplus_{1\le i\le n}{\Bbb R} \epsilon_i$ 
be the real Euclidean space with inner product 
$\langle\, ,\,\rangle$ such that
$\langle\epsilon_i,\epsilon_j\rangle=\delta_{ij}\,.$
Let 
$\Delta=\{\,
\pm\epsilon_i\pm\epsilon_j\,(\, 1\le i< j\le n\,), 
\pm 2\epsilon_i\,
(\, 1\le i\le n\,)\,\}$ be the root system of type $C_n,$
$\Delta^+=\{\,\epsilon_i\pm\epsilon_j\,
(\, 1\le i< j\le n\,), 2\epsilon_i\,
(\, 1\le i\le n\,)\,\}$ the set 
of positive roots, 
$\Pi=\{\,\alpha_i=\epsilon_i-\epsilon_{i+1}\,
(\, 1\le i\le n-1\,),\, \alpha_n=2\epsilon_n \,\}$  
the set of simple roots,
$P=\oplus_{1\le i\le n}~{\Bbb Z} \epsilon_i$  
the weight lattice, and 
$P^{\vee}=\oplus_{1\le i\le n}~{\Bbb Z} \epsilon_i
+{\Bbb Z}(\frac{1}{2}\sum_{i=1}^n \epsilon_i)$ 
the dual weight lattice for the root system $\Delta.$ 
We frequently write $\alpha\in\Delta^{+}$ as 
$\alpha >0.$ \par\smallskip
An element of the group algebra $A={\Bbb C}[P]$ is 
denoted by $e^{\lambda}\,,$ as is customary. 
Then the Weyl group $W=W(C_n)=\langle\,s_1, s_2,
\ldots, s_n\,\rangle$ (where each $s_i$ is a 
standard generator corresponding to the simple root
$\alpha_i$)  acts on $A$ as
$w(e^{\lambda})=e^{w\lambda}\;(w\in W).$   The symbol
$s_\alpha$ denoting the reflections is defined by
$s_\alpha(x)=x-\langle x,\alpha\rangle\alpha^{\vee}\,,$ 
with $\alpha^{\vee}=
2\alpha/{\langle\,\alpha, \alpha\rangle}\,$
for $x\in E$ and $\alpha\in\Delta\,.$\par\smallskip

The set of affine roots associated with $\Delta$ is 
$\Tilde{\Delta}=
\{\,\alpha+m\delta\,;\,\alpha\in\Delta,\, m\in{\Bbb Z}\},$ 
where $\delta$ denotes the constant function $1$ on $E.$ 
The simple roots are $a_0=-\theta+\delta$ with 
the highest root $\theta=2\epsilon_1\,$ and
$a_i=\alpha_i\in\Delta\,$ for  $1\le i\le n\,.$
We use the symbol introduced above, 
$s_i\,(0\le i\le n)$ to also represent the generator 
for the corresponding affine Weyl group.  We note that 
$s_0=\tau(\theta^\vee)s_{\theta}=
\tau(\epsilon_1)s_{2\epsilon_1}\,,$ 
where $\tau(\mu)$ is a translation
by $\mu\,.$\smallskip 

Let us introduce $V$ as the left free $A-$module of rank 
$|W|=2^n\,n!$ with the free basis $h_w\,(w\in W)\,;$ 
each element $F$ of $V$ can be written uniquely as 
$F=\sum_{w\in W} f_w h_w \;(f_w\in A).$
Then, let $A^{\sim}$ be a completion of the quotient 
field of $A\,.$ We then have
$V^{\sim}=A^{\sim}\otimes_A V\,.$
The action $r_w$ of the Weyl group $W$ on $V^{\sim}$ 
is defined by the following:
\begin{equation*}
r_w(f h_y)=w(f)h_{wy}\quad\text{ for }\quad f\in A\quad
\text{ and }\quad w,y\in W\,. 
\end{equation*}
\par\medskip\noindent
Moreover, the action of the translation 
$\tau(\mu)\,(\mu\in P^{\vee})$ 
for a parameter $u\in E$ is given by
$$\tau(\mu)e^\lambda=
q^{-\langle\,\lambda,\mu\,\rangle}e^\lambda
\quad\text{for}\quad\lambda\in P\,,\quad
\tau(\mu)h_w=q^{\langle\,\mu,wu\,\rangle}h_w
\quad\text{for}\quad w\in\,W$$
and
$$r_{\tau(\mu)}(f h_w)=\tau(\mu)(f)
q^{\langle\,\mu,wu\,\rangle}h_{w}\quad
\text{ for }\quad
f\in A\quad\text{ and }\quad w\in W\,.$$
This is an evaluation representation 
for which $e^\delta$ is identified with 
$q\,.$ \smallskip

Hereafter the symbol $r_w$ is used also to represent 
the element $w$ from the extended affine Weyl group
$W_{P^\vee}=W\ltimes{P^\vee}\,$
(the semidirect product of $W$ and $P^\vee$).
Then 
$r_w\, \tau({\epsilon_i})=
\tau(w({\epsilon_i}))\, r_w\,.$
Note also that, if 
$w=v\tau(\lambda),\; v\in W,\; \lambda\in P^{\vee}\,,$
we have 
$w(\mu)=v\mu-\langle\,\lambda, \mu\,\rangle\,\delta\;$ 
for $\mu\in P\,.$\smallskip

For an affine root $\alpha+m\delta\;
(\alpha\in\Delta\,, m\in\,{\Bbb Z}),$ 
define the $R$-matrix $R_{\alpha+m\delta}$ 
as an element of 
$End_{A^{\sim}}(V^{\sim})$ by the formula
\begin{equation*}
R_{\alpha+m\delta}h_y=
\begin{cases}
a_{\alpha+m\delta} h_y+
q^{m\,\langle\,\alpha^\vee,\, yu\,\rangle}
b_{\alpha+m\delta} h_{s_\alpha y}\,, 
& \quad   y^{-1}(\alpha)>0\,,\\[2pt]
c_{\alpha+m\delta} h_y+
q^{m\,\langle\,\alpha^\vee,\, yu\,\rangle}
d_{\alpha+m\delta} h_{s_\alpha y}\,,
& \quad  y^{-1}(\alpha)<0 
\end{cases}
\end{equation*}
for $y\in W\,,$ where 
{\allowdisplaybreaks
\begin{alignat*}{2}
a_{\alpha+m\delta}&=
\frac{1-q^{m}e^\alpha}{1-t_{\alpha}q^{m}e^\alpha} \,,
&\qquad 
b_{\alpha+m\delta}&=
\frac{1-t_{\alpha}}{1-t_{\alpha}q^{m}e^\alpha}\,,\\[2pt]
c_{\alpha+m\delta}&=
\frac{t_{\alpha}(1-q^{m}e^\alpha)}{1-t_{\alpha}q^{m}e^\alpha}\,,
&\qquad
d_{\alpha+m\delta}&=
\frac{q^{m}e^\alpha(1-t_{\alpha})}{1-t_{\alpha}q^{m}e^\alpha}
\end{alignat*}}\par\noindent
and $\alpha\mapsto\,t_\alpha$ is a 
$W$-invariant function taking positive values;
there are two different $t_\alpha\,,$ which we may 
write as $t_1=t_{\pm\epsilon_i\pm\epsilon_j},
t_2=t_{\pm2\epsilon_j}\,.$ \\ 
It is seen that
\begin{align}
r_{w}R_\alpha&=R_{w(\alpha)}r_{w} \quad 
\text{for}\quad\alpha\in\tilde{\Delta}\;,
w\in W_{P^\vee}\,,\label{eq2.1}\\[2pt]
R_{\beta}&=R_{-\beta}^{-1}\quad\text{for}
\quad\beta\in\tilde{\Delta}\label{eq2.2}
\end{align}
and
\begin{equation}
\begin{cases}
&R_{\epsilon_i-\epsilon_{j}}R_{\epsilon_i-\epsilon_{k}}
R_{\epsilon_j-\epsilon_{k}}
=R_{\epsilon_j-\epsilon_{k}}R_{\epsilon_i-\epsilon_{k}}
R_{\epsilon_i-\epsilon_{j}}\,,
\quad 1\le  i< j < k \le n\,,\\[2pt]
&R_{\epsilon_i-\epsilon_{j}}R_{2\epsilon_i}
R_{\epsilon_i+\epsilon_{j}}R_{2\epsilon_{j}}
=R_{2\epsilon_{j}}R_{\epsilon_i+\epsilon_{j}}
R_{2\epsilon_i}R_{\epsilon_i-\epsilon_{j}}\,, 
\quad 1\le i< j\le n\,.\\
\end{cases}\label{eq2.3}
\end{equation}
The relations in (\ref{eq2.3}) constitute 
the Yang-Baxter equation associated with the 
root system of type $C_n.$\medskip

Then we can state the definition of the 
QKZ equation for the root system of type $C_n.$ 
\begin{df}
The QKZ equation for the root system 
$C_{n}$ with a parameter 
$u=(u_1,\ldots,u_n)\in {\Bbb R}^n$ is 
the following system of equations: 
\begin{gather*}
r_{\tau(\epsilon_i)}^{-1}F=
R_{\tau(\epsilon_i)}F\,,\qquad 1\le i\le n,\\
\intertext{and}
r_{\tau(\frac{1}{2}(\epsilon_1+\cdots+\epsilon_n))}^{-1}F
=R_{\tau(\frac{1}{2}(\epsilon_1+\cdots+\epsilon_n))}F\,,
\end{gather*}
for $F\in\,V^{\sim}\,,$ with
\begin{align*}
R_{\tau(\epsilon_i)}&=
R_{\epsilon_i-\epsilon_{i-1}+\delta}\;\cdots\,
R_{\epsilon_i-\epsilon_1+\delta}\,
R_{2\epsilon_i+\delta}\,
R_{\epsilon_1+\epsilon_{i}}\cdots
R_{\epsilon_{i-1}+\epsilon_{i}}\\[2pt]
&\times
R_{\epsilon_{i}+\epsilon_{i+1}}\cdots
R_{\epsilon_{i}+\epsilon_{n}}\,
R_{2\epsilon_{i}}\,
R_{\epsilon_{i}-\epsilon_{n}}\cdots\,
R_{\epsilon_{i}-\epsilon_{i+1}}\\
\intertext{for $1\le i\le n\,,$ and}
R_{\tau(\frac{1}{2}(\epsilon_1+\cdots+\epsilon_n))}&=
(R_{2\epsilon_{1}}R_{\epsilon_{1}+
\epsilon_{2}}R_{\epsilon_{1}+\epsilon_{3}}
\cdots\,R_{\epsilon_{1}+\epsilon_{n}})\\[2pt]
&\times (R_{2\epsilon_{2}}R_{\epsilon_{2}+\epsilon_{3}}
\cdots\,R_{\epsilon_{2}+\epsilon_{n}})
\times\cdots\times 
(R_{2\epsilon_{n-1}}R_{\epsilon_{n-1}+\epsilon_{n}})
R_{2\epsilon_{n}}\,.
\end{align*}
\end{df}\medskip
{\it Remark.} If we introduce the operators 
$L_{\mu}\,(\mu\in\,P^{\vee})$ and 
$P^u_{\mu}\in\,End_{A^{\sim}}
(V^{\sim})\,(\mu\in\,P^{\vee},\,
u\in E)$ defined by
\begin{equation*}
L_{\mu} (\sum f_w h_w) =\sum L_{\mu}(f_w)h_w 
\quad\text{with}\quad
L_{\mu} (e^\lambda) 
= q^{\langle\mu,\lambda\rangle}e^\lambda
\quad (\lambda\in\,P)
\end{equation*}
and
\begin{equation*}
P^u_{\mu} (h_w) 
=q^{\langle \mu,wu\rangle} h_w \,,
\end{equation*}
then the equation above can be rewritten as
\begin{align*}
L_{\epsilon_i}F&=P^u_{\epsilon_i}R_{\tau(\epsilon_i)}F\\
\intertext{and}
L_{\frac{1}{2}(\epsilon_1+\cdots+\epsilon_n)}F&=
P^u_{\frac{1}{2}(\epsilon_1+\cdots+\epsilon_n)}
R_{\tau(\frac{1}{2}(\epsilon_1+\cdots+\epsilon_n))}F\,.
\end{align*}

Fulfilment of the compatibility condition of 
the QKZ equation is guaranteed by the Yang-Baxter 
equation (\ref{eq2.3}).\par\medskip
In the next section, we will construct a solution 
of the QKZ equation for the special case 
$u=-\lambda\epsilon_1\,(\lambda>0)$ through 
application of the $q$-Jordan-Pochhammer integral.
\section{Integrals and main result}\label{sec:3}
We introduce the form 
\begin{equation}
\Phi=x^{\lambda}\prod_{1\le j\le n}
\frac{ (ty_j/x)_\infty (ty_j^{-1}/x)_\infty}
{(y_j/x)_\infty (y_j^{-1}/x)_\infty}\frac{dx}{x}\,,
\end{equation}
where $(a)_\infty=\prod_{s\ge 0}(1-aq^s)\,.$
This can be regarded as a form of a restriction of 
the $q$-Jordan-Pochhammer integral 
$$x^{\lambda}\prod_{1\le j\le 2n}
\frac{ (ty_j/x)_\infty }
{(y_j/x)_\infty }\frac{dx}{x}\,,$$
which is studied in \cite{Mi1} and \cite{AKM}.

Next, to construct our solution in case of 
$u=-\lambda\epsilon_1\;(\lambda>0),$ 
we use the induced representaion of the 
Weyl group $W=W(C_n)$ from the trivial 
representation of a parabolic subgroup.

As a parabolic subgroup of $W,$ we choose a stabilizer
$W_{\epsilon_1}=\langle s_2, \ldots, s_n\rangle$ of 
$\epsilon_1.$  
A representative of the quotient
$W/W_{\epsilon_1}$ is fixed to be  
\begin{equation*}
\begin{cases}
&w_1=e,\, w_2=s_1,\, w_3=s_2s_1,\, \ldots,\, 
w_{n+1}=s_n\cdots\,s_2s_1,\\ 
&w_{n+2}=s_{n-1}w_{n+1},\,
w_{n+3}=s_{n-2}s_{n-1}w_{n+1},\,
\ldots,\, w_{2n}=s_{1}\cdots\, s_{n-1}w_{n+1}\,.
\end{cases}
\end{equation*} 

It is seen that the element 
$\overline{h}_e=\sum_{g\in\,W_{\epsilon_1}}h_g$
is invariant under the action of $W_{\epsilon_1}$ and the
induced representation of $W$ from ${\Bbb c}\overline{h}_e$
is given by the elements
 $$\overline{h}_{w_i}=
\sum_{g\in\,W_{\epsilon_1}}h_{w_i\,g}\quad(1\le\,i\le\,2n)\,.$$ 
 
Using $w_i$ as suffices, we define the 
following rational functions:

\begin{equation*}
\varphi_{w_i}=
\begin{cases}
\quad\frac{ \displaystyle\prod_{1\le \mu <i}
\biggl(1-\frac{y_{\mu}^{-1}}{x}\biggr)}
{ \displaystyle\prod_{1\le \mu \le i}
\biggl(1-t\frac{y_{\mu}^{-1}}{x}\biggr)}\,, 
&\quad 1\le i\le n\,,\\
\quad\frac{ \displaystyle\prod_{2n-i+1< \mu \le n}
\biggl(1-\frac{y_{\mu}}{x}\biggr)}
{ \displaystyle\prod_{2n-i+1\le \mu \le n}
\biggl(1-t\frac{y_{\mu}}{x}\biggr)}
{\displaystyle \prod_{1\le \mu \le n}}
\frac{ \displaystyle
\biggl(1-\frac{y_{\mu}^{-1}}{x}\biggr)}
{ \displaystyle
\biggl(1-t\frac{y_{\mu}^{-1}}{x}\biggr)}\,,
&\quad n+1\le i\le 2n\,.
\end{cases}
\end{equation*}
\medskip\noindent
Associated with the function $\Phi\,,$ we write 
\begin{equation*}
\langle \psi\rangle=\int_{\cal C}\psi\Phi 
\end{equation*}
for a rational function $\psi$ 
and a fixed cycle $\cal C\,,$ and define
the element $\Psi$ by
$$\Psi=\sum_{1\le i\le 2n}\langle\varphi_{w_i}\rangle
\overline{h}_{w_i}\,.$$
\medskip
Then we obtain the following, which will 
be proven in the next section.
\begin{prop}\label{prop3.1}
$r_{a_i}\Psi=R_{a_i}\Psi\qquad\text{for}
\qquad 0\le i\le n.$
\end{prop}

We are now in a position to state our main result.
\begin{thm}\label{thm3.2} The function
$$\Psi=\sum_{1\le i\le 2n}\langle\varphi_{w_i}\rangle
\overline{h}_{w_i}$$
satisfies the QKZ equation of type $C_n$ with
the parameter 
$u=-\lambda\epsilon_1 \,(\lambda>0)$ and $t_1=t_2=t\,:$
\begin{align}
r_{\tau(\epsilon_i)}^{-1}\Psi&=
R_{\tau(\epsilon_i)}\Psi\,,\qquad 1\le i\le n\,,\\
\intertext{and}
r_{\tau(\frac{1}{2}(\epsilon_1+\cdots+\epsilon_n))}^{-1}\Psi
&=R_{\tau(\frac{1}{2}(\epsilon_1+\cdots+\epsilon_n))}\Psi\,.
\end{align}
From this point we use the identification 
$y_i=e^{\epsilon_i}\,$ for $1\le i\le n\,.$
\end{thm}\par\medskip

It is seen that a system of fundamental solutions is 
obtained by taking suitable linearly 
independent cycles.\par\medskip
{\it Proof.} We first note 
$$r_{\tau(\epsilon_1)}^{-1}\Psi=r_{s_\theta\,s_0}\Psi\,.$$
Proposition 3.1  and (2.1) 
imply
$$r_{s_\theta s_0}\,\Psi
=r_{s_\theta}r_{s_0}\,\Psi
=r_{s_\theta}R_{\alpha_0}\,\Psi
=R_{s_\theta(\alpha_0)}r_{s_\theta}\,\Psi\,.
$$
Applying this process repeatedly, we finally obtain
\begin{align*}
r_{s_\theta s_0}\,\Psi
=&R_{s_\theta(\alpha_0)}R_{(s_1\cdots s_n)(s_{n-1}\cdots s_2)(\alpha_1)}
R_{(s_1\cdots s_n)(s_{n-1}\cdots s_3)(\alpha_2)}\cdots
R_{(s_1\cdots s_n)(\alpha_{n-1})}\\[2pt]
&\times
R_{(s_1\cdots s_{n-1})(\alpha_{n})}\cdots
R_{s_1(\alpha_{2})}R_{\alpha_{1}}\,\Psi\\[3pt]
=&R_{2\epsilon_1+\delta}\; R_{\epsilon_1+\epsilon_2}
\cdots\,R_{\epsilon_1+\epsilon_n}\;\,
R_{2\epsilon_1}\;\,R_{\epsilon_1-\epsilon_n}
\cdots\,R_{\epsilon_1-\epsilon_2}\,\Psi\,,
\end{align*}
since $ s_\theta=(s_1\cdots\,s_{n-1})(s_n\cdots\,s_1)\,.$ 
Thus we have
\begin{equation}
r_{\tau(\epsilon_1)}^{-1}\Psi
=R_{2\epsilon_1+\delta}\; R_{\epsilon_1+\epsilon_2}
\cdots\,R_{\epsilon_1+\epsilon_n}\;
R_{2\epsilon_1}\;R_{\epsilon_1-\epsilon_n}
\cdots\,R_{\epsilon_1-\epsilon_2}\,\Psi\,.\label{eq3.4}
\end{equation}
Next, let us apply $r_{s_{i-1}\cdots\,s_1}$ on 
both sides of (\ref{eq3.4}).
Then the left-hand side is
\begin{align*}
&r_{s_{i-1}\cdots\,s_1}r_{\,\tau(\epsilon_1)\,}^{-1}\Psi
=r_{\tau(\,s_{i-1}\cdots\,s_1(\epsilon_1)\,)}^{-1}
\,r_{s_{i-1}\cdots\,s_1}\Psi\\[2pt]
&=r_{\tau(\,s_{i-1}\cdots\,s_1(\epsilon_1)\,)}^{-1}\,
R_{s_{i-1}\cdots\,s_2(\alpha_1)}\,R_{s_{i-1}\cdots\,s_3(\alpha_2)}\,
\cdots\,
R_{s_{i-1}(\alpha_{i-2})}R_{\alpha_{i-1}}\Psi\\[2pt]
&=r_{\tau(\epsilon_i)}^{-1}
R_{\epsilon_1-\epsilon_i}R_{\epsilon_2-\epsilon_i}\,
\cdots\,
R_{\epsilon_{i-2}-\epsilon_i}R_{\epsilon_{i-1}-\epsilon_i}\Psi\\[2pt]
&=R_{\epsilon_1-\epsilon_i-\delta}\,R_{\epsilon_2-\epsilon_i-\delta}\,
\cdots\,
R_{\epsilon_{i-2}-\epsilon_i-\delta}\,
R_{\epsilon_{i-1}-\epsilon_i-\delta}\;
r_{\tau(\epsilon_i)}^{-1}\Psi\,.
\end{align*}
This follows from the relation 
$\tau(-\epsilon_i)(\epsilon_j-\epsilon_i)
=\epsilon_j-\epsilon_i-\delta\,.$\par\smallskip
On the other hand, the right-hand side is
\begin{align*}
&r_{s_{i-1}\cdots\,s_1}
R_{2\epsilon_1+\delta}\, R_{\epsilon_1+\epsilon_2}
\cdots\,R_{\epsilon_1+\epsilon_n}\,
R_{2\epsilon_1}\,R_{\epsilon_1-\epsilon_n}
\cdots\,R_{\epsilon_1-\epsilon_2}\,\Psi\\[2pt]
=&R_{2\epsilon_i+\delta}\, R_{\epsilon_1+\epsilon_i}
R_{\epsilon_2+\epsilon_i}\cdots\,R_{\epsilon_{i-1}+\epsilon_i}
\; R_{\epsilon_{i}+\epsilon_{i+1}}\cdots
R_{\epsilon_{i}+\epsilon_{n}}\,
R_{2\epsilon_i}\,\Psi\,.
\end{align*}
Here we have used
$$r_{s_{i-1}\cdots\,s_1}^{-1}\,\Psi=
R_{\epsilon_{i}-\epsilon_{n}}\cdots 
R_{\epsilon_{n-1}-\epsilon_{n}}\,\Psi\,.$$
Therefore we reach the desired relation (3.2) 
by using (\ref{eq2.2}).\par\medskip

Next we proceed to derive (3.3).\medskip

For $1\le i\le n\,,$ we have
{\allowdisplaybreaks
\begin{equation}
\begin{split}
&r_{\tau(\frac{1}{2}(\epsilon_1+\cdots+\epsilon_n))}^{-1}
\langle\,\varphi_{w_i}\,\rangle \\
&=\int_C\,x^{\lambda}\prod_{k=1}^n
\frac{\displaystyle\biggl(q^{\frac{1}{2}}\frac{ty_k}{x}\biggr)_\infty}
{\displaystyle\biggl(q^{\frac{1}{2}}\frac{y_k}{x}\biggr)_\infty}
\frac{\displaystyle\;\;
\prod_{k=i+1}^n\,\biggl(q^{-\frac{1}{2}}\,\frac{ty_k^{-1}}{x}\biggr)_\infty\,
\prod_{k=1}^i\,\biggl(q^{\frac{1}{2}}\,\frac{ty_k^{-1}}{x}\biggr)_\infty}
{\displaystyle\quad
\prod_{k=i}^n\;\;\biggl(q^{-\frac{1}{2}}\,\frac{y_k^{-1}}{x}\biggr)_\infty\;\;
\prod_{k=1}^{i-1}\,\biggl(q^{\frac{1}{2}}\,\frac{y_k^{-1}}{x}\biggr)_\infty}
\frac{dx}{x}
\end{split}
\end{equation}}
and
{\allowdisplaybreaks
\begin{equation}
\begin{split}
&r_{\tau(\frac{1}{2}(\epsilon_1+\cdots+\epsilon_n))}^{-1}
\langle\,\varphi_{w_{n+i}}\,\rangle \\
&=\int_C\,x^{\lambda}\;
\frac{\displaystyle\;\;
\prod_{k=1}^{n-i}\,\biggl(q^{\frac{1}{2}}\,\frac{ty_k}{x}\biggr)_\infty\;
\prod_{k=n-i+1}^n\,\biggl(q^{\frac{3}{2}}\,\frac{ty_k}{x}\biggr)_\infty}
{\displaystyle
\prod_{k=1}^{n-i+1}\,\biggl(q^{\frac{1}{2}}\,\frac{y_k}{x}\biggr)_\infty\;
\prod_{k=n-i+2}^{n}\,\biggl(q^{\frac{3}{2}}\,\frac{y_k}{x}\biggr)_\infty}
\prod_{k=1}^n
\frac{\displaystyle\biggl(q^{\frac{1}{2}}\frac{ty_k^{-1}}{x}\biggr)_\infty}
{\displaystyle\biggl(q^{\frac{1}{2}}\frac{y_k^{-1}}{x}\biggr)_\infty}
\frac{dx}{x}\,.
\end{split}
\end{equation}}
By changing the integration variable such that 
$x\mapsto q^{-1/2}x$, from (3.5) we have
{\allowdisplaybreaks
\begin{equation*}
\begin{split}
&r_{\tau(\frac{1}{2}(\epsilon_1+\cdots+\epsilon_n))}^{-1}
\langle\,\varphi_{w_i}\,\rangle \\
&=q^{-\frac{\lambda}{2}}\int_C\,x^{\lambda}\prod_{k=1}^n
\frac{\displaystyle\biggl(q\frac{ty_k}{x}\biggr)_\infty}
{\displaystyle\biggl(q\frac{y_k}{x}\;\biggr)_\infty}
\frac{\displaystyle\;\;
\prod_{k=i+1}^n\,\biggl(\,\frac{ty_k^{-1}}{x}\biggr)_\infty\,
\prod_{k=1}^i\,\biggl(q\,\frac{ty_k^{-1}}{x}\biggr)_\infty}
{\displaystyle\quad
\prod_{k=i}^n\;\;\biggl(\,\frac{y_k^{-1}}{x}\,\biggr)_\infty\;\;
\prod_{k=1}^{i-1}\,\biggl(q\,\frac{y_k^{-1}}{x}\,\biggr)_\infty}
\frac{dx}{x}\\
&=q^{-\frac{\lambda}{2}}\langle\;g\varphi_{w_{n+i}}\;\rangle
\end{split}
\end{equation*}}\noindent
with
$$g=s_n(s_{n-1}s_n)(s_{n-2}s_{n-1}s_n)\cdots(s_1\cdots\,s_n)\in W\,.$$ 
Here we note $g(\epsilon_i)=-\epsilon_{n-i+1}$ for each 
$1\le i\le n\,.$\par\bigskip

Similarly, as a result of the change $x\mapsto q^{1/2}x$, 
from (3.6) we have
{\allowdisplaybreaks
\begin{equation*}
\begin{split}
&r_{\tau(\frac{1}{2}(\epsilon_1+\cdots+\epsilon_n))}^{-1}
\langle\,\varphi_{w_{n+i}}\,\rangle \\
&=q^{\frac{\lambda}{2}}\int_C\,x^{\lambda}\;
\frac{\displaystyle\;\;
\prod_{k=1}^{n-i}\;\biggl(\frac{ty_k}{x}\biggr)_\infty\;\,
\prod_{k=n-i+1}^n\,\biggl(q\,\frac{ty_k}{x}\biggr)_\infty}
{\displaystyle
\prod_{k=1}^{n-i+1}\,\biggl(\,\frac{y_k}{x}\biggr)_\infty\;
\prod_{k=n-i+2}^{n}\,\biggl(q\,\frac{y_k}{x}\biggr)_\infty}
\prod_{k=1}^n
\frac{\displaystyle\biggl(\frac{ty_k^{-1}}{x}\biggr)_\infty}
{\displaystyle\biggl(\frac{y_k^{-1}}{x}\biggr)_\infty}
\frac{dx}{x}\\
&=q^{\frac{\lambda}{2}}\langle\;g\varphi_{w_{i}}\;\rangle
\end{split}
\end{equation*}}
with the same $g\in W\,.$

As for this 
$g=s_n(s_{n-1}s_n)(s_{n-2}s_{n-1}s_n)\cdots(s_1\cdots\,s_n)\in W\,,$ 
we have
\begin{align*}
&gw_i=w_{n+i}\,s_n(s_{n-1}s_n)(s_{n-2}s_{n-1}s_n)
\cdots (s_{2}\cdots\ s_{n-1}s_n)\,,\\[2pt]
&gw_{n+i}=w_{i}\,s_n(s_{n-1}s_n)(s_{n-2}s_{n-1}s_n)
\cdots (s_{2}\cdots\ s_{n-1}s_n)
\end{align*}\noindent
for $1\le i\le n\,.$ These relations lead to 
\begin{align*}
&g\overline{h}_{w_i}=\overline{h}_{gw_i}
=\overline{h}_{w_{n+i}}\,,\\[2pt]
&g\overline{h}_{w_{n+i}}=\overline{h}_{gw_{n+i}}
=\overline{h}_{w_{i}}\,
\end{align*}
for $1\le i\le n\,.$\medskip

On the other hand, noting $u=-\lambda\epsilon_1\,,$ we obtain
{\allowdisplaybreaks
\begin{align*}
&\tau(-\frac{1}{2}
(\epsilon_1+\cdots+\epsilon_n))\overline{h}_{w_i}
=q^{\langle\,-\frac{1}{2}(\epsilon_1+\cdots+\epsilon_n),\,
-\lambda\epsilon_1\,\rangle} \overline{h}_{w_i}
=q^{\frac{\lambda}{2}}\overline{h}_{w_i}\,,\\[2pt]
&\tau(-\frac{1}{2}(\epsilon_1+\cdots+\epsilon_n))\overline{h}_{w_{n+i}}
=q^{\langle\,-\frac{1}{2}(\epsilon_1+\cdots+\epsilon_n),\,
\lambda\epsilon_{n-i+1}\,\rangle} \overline{h}_{w_{n+i}}=
q^{-\frac{\lambda}{2}}\overline{h}_{w_{n+i}}
\end{align*}}
for $1\le i\le n\,.$\medskip

Combining these relations, we get
{\allowdisplaybreaks
\begin{align*}
&\quad\tau(-\frac{1}{2}(\epsilon_1+\cdots+\epsilon_n))\,\Psi\\[2pt]
&=\tau(-\frac{1}{2}(\epsilon_1+\cdots+\epsilon_n))
\sum_{1\le i\le n}
\biggl\{\langle\, \varphi_{w_i}\,\rangle\,\overline{h}_{w_{i}}
+\langle\, \varphi_{w_{n+i}}\,\rangle\,\overline{h}_{w_{n+i}}\biggr\}\\[2pt]
&=\sum_{1\le i\le n}\biggl\{
q^{-\frac{\lambda}{2}}\,\langle\, g\varphi_{w_{n+i}}\,\rangle\,
q^{\frac{\lambda}{2}}\overline{h}_{w_{i}}+
q^{\frac{\lambda}{2}}\,\langle\, g\varphi_{w_{i}}\,\rangle\,
q^{-\frac{\lambda}{2}}\overline{h}_{w_{n+i}}\biggr\}\\[2pt]
&=\sum_{1\le i\le n}\biggl\{
\,\langle\, g\varphi_{w_{n+i}}\,\rangle\,
\overline{h}_{w_{i}}+
\,\langle\, g\varphi_{w_{i}}\,\rangle\,
\overline{h}_{w_{n+i}}\biggr\}\\[2pt]
&=\sum_{1\le i\le n}\biggl\{
\,\langle\, g\varphi_{w_{n+i}}\,\rangle\,
g\overline{h}_{w_{n+i}}+
\,\langle\, g\varphi_{w_{i}}\,\rangle\,
g\overline{h}_{w_{i}}\biggr\}\\[2pt]
&=r_g\,\Psi\,.
\end{align*}}
\medskip\noindent
At this stage, applying the relation
{\allowdisplaybreaks
\begin{align*}
&r_{(s_n(s_{n-1}s_n)\cdots(s_{k+1}\cdots s_{n}))s_k\cdots s_n}\,\Psi\\[2pt]
&=R_{(s_n(s_{n-1}s_n)\cdots(s_{k+1}\cdots s_{n}))s_k\cdots s_{n-1}(\alpha_n)}
R_{(s_n(s_{n-1}s_n)\cdots(s_{k+1}\cdots s_{n}))s_k\cdots s_{n-2}(\alpha_{n-1})}\\[2pt]
&\times\cdots
\times R_{(s_n(s_{n-1}s_n)\cdots(s_{k+1}\cdots s_{n}))s_k(\alpha_{k+1})}
R_{(s_n(s_{n-1}s_n)\cdots(s_{k+1}\cdots s_{n}))(\alpha_{k})}\\
&\times r_{(s_n(s_{n-1}s_n)\cdots(s_{k+1}\cdots s_{n-1}s_n))}\,\Psi\\[3pt]
&=R_{2\epsilon_k}R_{\epsilon_{k}+\epsilon_{k+1}}\cdots
R_{\epsilon_{k}+\epsilon_{n-1}}R_{\epsilon_{k}+\epsilon_{n}}\\[2pt]
&\times r_{(s_n(s_{n-1}s_n)\cdots(s_{k+1}\cdots s_{n-1}s_n))}\,\Psi\,,\qquad (1\le k\le n)
\end{align*}}
repeatedly, we finally obtain 
{\allowdisplaybreaks
\begin{align*}
r_{g}\,\Psi
&=(R_{\,2\epsilon_1}R_{\,\epsilon_{1}+\epsilon_{2}}\cdots
R_{\,\epsilon_{1}+\epsilon_{n}})
(R_{\,2\epsilon_2}R_{\,\epsilon_{2}+\epsilon_{3}}\cdots
R_{\,\epsilon_{2}+\epsilon_{n}})\\[2pt]
&\times\cdots
\times (R_{\,2\epsilon_{n-1}\,}R_{\,\epsilon_{n-1}+\epsilon_{n}\,})\,
R_{\,2\epsilon_{n-1}}\,\Psi\,.
\end{align*}}
Therefore, we reach the desired result (3.3).\qquad\qed\bigskip
\section{Proof of Proposition 3.1 }\label{sec4}

To prove Proposition 3.1, we start by considering
the action of $s_i\in W$ on the $\varphi_{w_k}\,.$ 
\begin{lem}\label{lem4.1}
\begin{enumerate}
\renewcommand{\labelenumi}{(\alph{enumi})}
\item If $1\le i \le n-1\,,\; s_i\,\varphi_{w_k}=\varphi_{w_k}$ 
for each \;$1\le\,k\le 2n$ such that 
$k\ne i,\,i+1,\,2n-i,\,2n-i+1\,.$\smallskip 
\item $s_n\,\varphi_{w_k}=\varphi_{w_k}$ 
for each \;$1\le\,k\le 2n$ such that $k\ne n,\,n+1\,.$ \smallskip
\item $s_0\,\varphi_{w_k}=\varphi_{w_k}$ 
for each \;$1\le\,k\le 2n$ such that $k\ne 1,\,2n\,.$ 
\end{enumerate}
\end{lem}
{\it Proof.} These assertions follow from the definition 
of $s_i$ and $\varphi_{w_k}\,.$\qquad\qed\par\medskip
Moreover we have 
\begin{lem}\label{lem4.2}
\begin{enumerate}
\renewcommand{\labelenumi}{(\alph{enumi})}
\item For $1\le i\le n-1\,;$  
\begin{equation}
\begin{cases}
\quad s_i\varphi_{w_{i+1}}&=
a_{\alpha_i}\;\varphi_{w_i}+d_{\alpha_i}\;\varphi_{w_{i+1}}\,,\\
\quad s_i\varphi_{w_i}&=
b_{\alpha_i}\;\varphi_{w_i}+c_{\alpha_i}\;\varphi_{w_{i+1}}\,,\\
\end{cases}
\end{equation}\noindent
and
\begin{equation}
\begin{cases}
\quad s_i\varphi_{w_{2n-i+1}}&=
a_{\alpha_i}\;\varphi_{w_{2n-i}}+d_{\alpha_i}\;\varphi_{w_{2n-i+1}}\,,\\
\quad s_i\varphi_{w_{2n-i}}&=
b_{\alpha_i}\;\varphi_{w_{2n-i}}+c_{\alpha_i}\;\varphi_{w_{2n-i+1}}\,.\\
\end{cases}
\end{equation}
\par\medskip\noindent
\item
\begin{equation}
\begin{cases}
\quad s_n\varphi_{w_{n+1}}&=
a_{\alpha_n}\;\varphi_{w_n}+d_{\alpha_n}\;\varphi_{w_{n+1}}\,,\\
\quad s_n\varphi_{w_n}&=
b_{\alpha_n}\;\varphi_{w_n}+c_{\alpha_n}\;\varphi_{w_{n+1}}\,.\\
\end{cases}
\end{equation}
\end{enumerate}
\end{lem}

{\it Proof.} By direct calculation or 
expansion of partial fractions, 
we find 
{\allowdisplaybreaks
\begin{equation}
\begin{split}
&\frac{ \displaystyle 1-\frac{y_{i+1}^{-1}}{x}}
{\displaystyle 
\biggl(1-t\frac{y_{i}^{-1}}{x}\biggr)
\biggl(1-t\frac{y_{i+1}^{-1}}{x}\biggr)}\\
&=
a_{\alpha_i}\,\frac{1}{\displaystyle 
1-t\frac{y_{i}^{-1}}{x}}
+d_{\alpha_{i}}\,
\frac{ \displaystyle 
1-\frac{y_{i}^{-1}}{x}}
{ \displaystyle
\biggl(1-t\frac{y_{i+1}^{-1}}{x}\biggr)
\biggl(1-t\frac{y_{i}^{-1}}{x}\biggr)}
\end{split}
\end{equation}}
and
{\allowdisplaybreaks
\begin{equation}
\frac{1}{\displaystyle 1-t\frac{y_{i+1}^{-1}}{x}}
=b_{\alpha_i}\,\frac{1}{\displaystyle 1-t\frac{y_{i}^{-1}}{x}}
+c_{\alpha_i}\,
\frac{ \displaystyle 1-\frac{y_{i}^{-1}}{x}}
{ \displaystyle
\biggl(1-t\frac{y_{i+1}^{-1}}{x}\biggr)
\biggl(1-t\frac{y_{i}^{-1}}{x}\biggr)}\,.
\end{equation}}
Multiplying the factor
$$\prod_{j=1}^{i-1}
\frac{\displaystyle 1-\;\frac{y_{j}^{-1}}{x}}
{\displaystyle 1-t\frac{y_{j}^{-1}}{x}}$$
on both sides of each equality (4.4) or (4.5), 
we get the desired relations (4.1).

While the change of variables 
$\epsilon_i\mapsto-\epsilon_{i+1}$ and
$\epsilon_{i+1}\mapsto-\epsilon_{i}$ leave 
$\alpha_i$ unchanged, they produce the 
following from (4.4) and (4.5):
{\allowdisplaybreaks
\begin{equation}
\begin{split}
&\frac{\displaystyle 1-\frac{y_{i}}{x}}
{\displaystyle 
\biggl(1-t\frac{y_{i}}{x}\biggr)
\biggl(1-t\frac{y_{i+1}}{x}\biggr)}\\
&=
a_{\alpha_i}\,\frac{1}{\displaystyle 
1-t\frac{y_{i+1}}{x}}
+d_{\alpha_{i}}\,
\frac{ \displaystyle 
1-\frac{y_{i+1}}{x}}
{ \displaystyle
\biggl(1-t\frac{y_{i}}{x}\biggr)
\biggl(1-t\frac{y_{i+1}}{x}\biggr)}
\end{split}
\end{equation}}
and
{\allowdisplaybreaks
\begin{equation}
\frac{1}{\displaystyle 1-t\frac{y_{i}}{x}}
=b_{\alpha_i}\,\frac{1}{\displaystyle 1-t\frac{y_{i+1}}{x}}
+c_{\alpha_i}\,
\frac{ \displaystyle 1-\frac{y_{i+1}}{x}}
{ \displaystyle
\biggl(1-t\frac{y_{i}}{x}\biggr)
\biggl(1-t\frac{y_{i+1}}{x}\biggr)}\,.
\end{equation}} 

Multiplying the factor
$$
\prod_{j=i+2}^{n}
\frac{\displaystyle 1-\;\frac{y_{j}}{x}}{\displaystyle 1-t\frac{y_{j}}{x}}
\prod_{j=1}^{n}
\frac{\displaystyle 1-\;\frac{y_{j}^{-1}}{x}}
{\displaystyle 1-t\frac{y_{j}^{-1}}{x}}$$
on both sides of equalities (4.6) and (4.7), 
we obtain the desired relations (4.2).
 
Similarly, changing $\epsilon_i\mapsto\epsilon_{n}$ and
$\epsilon_{i+1}\mapsto-\epsilon_{n}$ induces 
$\alpha_i\mapsto\alpha_n $ and 
leads from equalities (4.4) and (4.5) to
{\allowdisplaybreaks
\begin{equation}
\begin{split}
&\frac{\displaystyle 1-\frac{y_{n}}{x}}
{\displaystyle 
\biggl(1-t\frac{y_{n}}{x}\biggr)
\biggl(1-t\frac{y_{n}^{-1}}{x}\biggr)}\\
&=
a_{\alpha_n}\,\frac{1}{\displaystyle 
1-t\frac{y_{n}^{-1}}{x}}
+d_{\alpha_{n}}\,
\frac{ \displaystyle 
1-\frac{y_{n}^{-1}}{x}}
{ \displaystyle
\biggl(1-t\frac{y_{n}}{x}\biggr)
\biggl(1-t\frac{y_{n}^{-1}}{x}\biggr)}
\end{split}
\end{equation}}
and
{\allowdisplaybreaks
\begin{equation}
\frac{1}{\displaystyle 1-t\frac{y_{n}}{x}}
=b_{\alpha_n}\,\frac{1}{\displaystyle 1-t\frac{y_{n}^{-1}}{x}}
+c_{\alpha_n}\,
\frac{ \displaystyle 1-\frac{y_{n}^{-1}}{x}}
{ \displaystyle
\biggl(1-t\frac{y_{n}}{x}\biggr)
\biggl(1-t\frac{y_{n}^{-1}}{x}\biggr)}\,.
\end{equation}} 
Multiplying the factor
$$
\prod_{j=1}^{n-1}
\frac{\displaystyle 1-\;\frac{y_{j}^{-1}}{x}}
{\displaystyle 1-t\frac{y_{j}^{-1}}{x}}$$
on both sides of equalities (4.8) and (4.9), 
we get the desired relations (4.3).\quad$\square$\par\medskip
In contrast to the action of $s_i$ for $1\le i\le n\,,$ the action of
$s_0$ is understood as it acts on the $q$ de Rham cohomology, 
not on the rational functions.
\begin{lem}\label{lem4.3}
\begin{equation}
\begin{cases}
\quad q^{\lambda}\langle\; s_0\varphi_{w_{1}}\;\rangle &=
a_{\delta-\theta}\langle\;\varphi_{w_{2n}}\;\rangle+
q^{\lambda}d_{\delta-\theta}\langle\;\varphi_{w_{1}}\;\rangle\,,\\
\quad q^{-\lambda}\langle\; s_0\varphi_{w_{2n}}\;\rangle&=
q^{-\lambda}b_{\delta-\theta}\langle\;\varphi_{w_{2n}}\;\rangle
+c_{\delta-\theta}\langle\;\varphi_{w_{1}}\;\rangle\,.\\
\end{cases}
\end{equation}
\end{lem}
{\it Proof.}
Make the change of variables $\epsilon_i\mapsto -\epsilon_{1}$ and
$\epsilon_{i+1}\mapsto \epsilon_{1}-\delta $ 
(i.e. $y_i^{-1}\mapsto y_1\,, \;y_{i+1}^{-1}\mapsto qy_1^{-1}$)
in (4.4) and (4.5).
Then we have
{\allowdisplaybreaks
\begin{equation}
\begin{split}
&\frac{\displaystyle 1-q\frac{y_{1}^{-1}}{x}}
{\displaystyle 
\biggl(1-t\frac{y_{1}}{x}\biggr)
\biggl(1-tq\frac{y_{1}^{-1}}{x}\biggr)}\\
&=
a_{\delta-\theta}\,\frac{1}{\displaystyle 
1-t\frac{y_{1}}{x}}
+d_{\delta-\theta}\,
\frac{ \displaystyle 
1-\frac{y_{1}}{x}}
{ \displaystyle
\biggl(1-tq\frac{y_{1}^{-1}}{x}\biggr)
\biggl(1-t\frac{y_{1}}{x}\biggr)}
\end{split}
\end{equation}}
and
{\allowdisplaybreaks
\begin{equation}
\frac{1}{\displaystyle 1-tq\frac{y_{1}^{-1}}{x}}
=b_{\delta-\theta}\,\frac{1}{\displaystyle 1-t\frac{y_{1}}{x}}
+c_{\delta-\theta}\,
\frac{ \displaystyle 1-\frac{y_{1}}{x}}
{ \displaystyle
\biggl(1-tq\frac{y_{1}^{-1}}{x}\biggr)
\biggl(1-t\frac{y_{1}}{x}\biggr)}\,.
\end{equation}} 
Integration after
multiplying the factor
$$
\prod_{j=2}^{n}
\frac{\displaystyle 1-\;\frac{y_{j}}{x}}
{\displaystyle 1-t\frac{y_{j}}{x}}
\prod_{j=1}^{n}
\frac{\displaystyle 1-\;\frac{y_{j}^{-1}}{x}}
{\displaystyle 1-t\frac{y_{j}^{-1}}{x}}\;\Phi$$
on both sides of equalities (4.11) and (4.12) gives 
the following:
{\allowdisplaybreaks
\begin{align}
& \int_C\,x^{\lambda}
\frac{\displaystyle\prod_{k=1}^n\biggl(q\frac{ty_k}{x}\biggr)_\infty}
{\displaystyle\biggl(\frac{y_1}{x}\biggr)_\infty\prod_{k=2}^n
\biggl(q\frac{y_k}{x}\biggr)_\infty}
\frac{\displaystyle\biggl(q^2\,\frac{ty_1^{-1}}{x}\biggr)_\infty}
{\displaystyle\biggl(q^2\,\frac{y_1^{-1}}{x}\biggr)_\infty}
\prod_{k=2}^n
\frac{\displaystyle
\biggl(q\,\frac{ty_k^{-1}}{x}\biggr)_\infty}
{\displaystyle
\biggl(q\;\frac{y_k^{-1}}{x}\biggr)_\infty}
\frac{dx}{x} \\
&=a_{\delta-\theta}\,\int_C\,x^{\lambda}
\frac{\displaystyle\prod_{k=1}^n\biggl(q\frac{ty_k}{x}\biggr)_\infty}
{\displaystyle\biggl(\frac{y_1}{x}\biggr)_\infty\prod_{k=2}^n
\biggl(q\frac{y_k}{x}\biggr)_\infty}
\prod_{k=1}^n
\frac{\displaystyle
\biggl(q\,\frac{ty_k^{-1}}{x}\,\biggr)_\infty}
{\displaystyle
\biggl(q\;\frac{y_k^{-1}}{x}\,\biggr)_\infty}
\frac{dx}{x}\notag\\
&+d_{\delta-\theta}\,\int_C\,x^{\lambda}
\prod_{k=1}^n
\frac{\displaystyle
\biggl(q\,\frac{ty_k}{x}\biggr)_\infty}
{\displaystyle
\biggl(q\,\frac{y_k}{x}\,\biggr)_\infty}
\frac{\displaystyle
\biggl(q^2 \,\frac{y_1^{-1}}{x}\biggr)_\infty
\prod_{k=2}^n
\biggl(q\,\frac{ty_k^{-1}}{x}\biggr)_\infty}{\displaystyle
\prod_{k=1}^n\biggl(q\,\frac{y_k^{-1}}{x}\,\biggr)_\infty}
\frac{dx}{x}\notag
\end{align}}
and
{\allowdisplaybreaks
\begin{align}
&\int_C\,x^{\lambda}
\frac{\displaystyle\biggl(\frac{ty_1}{x}\biggr)_\infty}
{\displaystyle\biggl(\frac{y_1}{x}\biggr)_\infty}
\prod_{k=2}^n
\frac{\displaystyle
\biggl(q\,\frac{ty_k}{x}\biggr)_\infty}{\displaystyle
\biggl(q\,\frac{y_k}{x}\biggr)_\infty}
\frac{\displaystyle
\biggl(q^2 \frac{y_1^{-1}}{x}\,\biggr)_\infty
\prod_{k=2}^n
\biggl(q\frac{ty_k^{-1}}{x}\,\biggr)_\infty}{\displaystyle
\prod_{k=1}^n\biggl(q\frac{y_k^{-1}}{x}\,\biggr)_\infty}
\frac{dx}{x} \\
&=b_{\delta-\theta}\,\int_C\,x^{\lambda}
\frac{\displaystyle\prod_{k=1}^n\biggl(q\frac{ty_k}{x}\biggr)_\infty}
{\displaystyle\biggl(\frac{y_1}{x}\biggr)_\infty\prod_{k=2}^n
\biggl(q\frac{y_k}{x}\biggr)_\infty}
\prod_{k=1}^n
\frac{\displaystyle
\biggl(q\,\frac{ty_k^{-1}}{x}\,\biggr)_\infty}
{\displaystyle
\biggl(q\,\frac{y_k^{-1}}{x}\,\biggr)_\infty}
\frac{dx}{x}\notag\\
&+c_{\delta-\theta}\,\int_C\,x^{\lambda}
\prod_{k=1}^n
\frac{\displaystyle
\biggl(q\,\frac{ty_k}{x}\;\biggr)_\infty}
{\displaystyle
\biggl(q\,\frac{y_k}{x}\;\biggr)_\infty}
\frac{\displaystyle
\biggl(q^2 \frac{y_1^{-1}}{x}\biggr)_\infty
\prod_{k=2}^n
\biggl(q\frac{ty_k^{-1}}{x}\biggr)_\infty}{\displaystyle
\prod_{k=1}^n\biggl(q\frac{y_k^{-1}}{x}\biggr)_\infty}
\frac{dx}{x}\notag
\end{align}}
Here, changing the integration variable such that 
$x\mapsto qx\,,$ we have
{\allowdisplaybreaks
\begin{align*}
&\int_C\,x^{\lambda}
\frac{\displaystyle\prod_{k=1}^n\biggl(q\frac{ty_k}{x}\biggr)_\infty}
{\displaystyle\biggl(\frac{y_1}{x}\biggr)_\infty\prod_{k=2}^n
\biggl(q\frac{y_k}{x}\biggr)_\infty}
\frac{\displaystyle\biggl(q^2\,\frac{ty_1^{-1}}{x}\,\biggr)_\infty}
{\displaystyle\biggl(q^2\,\frac{y_1^{-1}}{x}\;\biggr)_\infty}
\prod_{k=2}^n
\frac{\displaystyle
\biggl(q\,\frac{ty_k^{-1}}{x}\,\biggr)_\infty}
{\displaystyle
\biggl(q\,\frac{y_k^{-1}}{x}\;\biggr)_\infty}
\frac{dx}{x}\\
&=q^\lambda\;\int_C\,x^{\lambda}
\frac{\displaystyle\prod_{k=1}^n\biggl(\frac{ty_k}{x}\biggr)_\infty}
{\displaystyle\biggl(q^{-1}\,\frac{y_1}{x}\biggr)_\infty\prod_{k=2}^n
\biggl(\frac{y_k}{x}\biggr)_\infty}
\frac{\displaystyle\biggl(q\,\frac{ty_1^{-1}}{x}\,\biggr)_\infty}
{\displaystyle\biggl(q\,\frac{y_1^{-1}}{x}\;\biggr)_\infty}
\prod_{k=2}^n
\frac{\displaystyle
\biggl(\,\frac{ty_k^{-1}}{x}\,\biggr)_\infty}
{\displaystyle
\biggl(\,\frac{y_k^{-1}}{x}\;\biggr)_\infty}
\frac{dx}{x}\\
&=q^\lambda \,\langle\, s_0\varphi_{w_1}\,\rangle
\end{align*}}
and
{\allowdisplaybreaks
\begin{align*}
&\int_C\,x^{\lambda}
\prod_{k=1}^n
\frac{\displaystyle
\biggl(q\,\frac{ty_k}{x}\biggr)_\infty}
{\displaystyle
\biggl(q\,\frac{y_k}{x}\biggr)_\infty}
\frac{\displaystyle
\biggl(q^2 \,\frac{y_1^{-1}}{x}\biggr)_\infty
\prod_{k=2}^n
\biggl(q\,\frac{ty_k^{-1}}{x}\biggr)_\infty}{\displaystyle
\prod_{k=1}^n\biggl(q\,\frac{y_k^{-1}}{x}\biggr)_\infty}
\frac{dx}{x}\\
&=q^\lambda\;\int_C\,x^{\lambda}
\prod_{k=1}^n
\frac{\displaystyle
\biggl(\,\frac{ty_k}{x}\biggr)_\infty}
{\displaystyle
\biggl(\,\frac{y_k}{x}\biggr)_\infty}
\frac{\displaystyle
\biggl(q \frac{y_1^{-1}}{x}\biggr)_\infty
\prod_{k=2}^n
\biggl(\frac{ty_k^{-1}}{x}\biggr)_\infty}{\displaystyle
\prod_{k=1}^n\biggl(\frac{y_k^{-1}}{x}\biggr)_\infty}
\frac{dx}{x}\\
&=q^\lambda \,\langle\, \varphi_{w_1}\,\rangle\,.
\end{align*}}
Therefore, it is seen that (4.13) and (4.14) are equivalent to the 
desired relations (4.10).
\qquad\qed\medskip

Next, we consider the asction of $W$ on the $\overline{h}_{w_k}\,.$
\begin{lem}\label{lem4.4}
\begin{enumerate}
\renewcommand{\labelenumi}{(\alph{enumi})}
\item If $1\le i\le n-1\,,$
$\overline{h}_{s_iw_{k}}=\overline{h}_{w_{k}}$
for $k\ne i,i+1,2n-i,2n-i+1\,.$
\item
$\overline{h}_{s_nw_{k}}=\overline{h}_{w_{k}}$
for $k\ne n, n+1\,.$
\item
$\overline{h}_{s_\theta w_{k}}=\overline{h}_{w_{k}}$
for $k\ne 1, 2n\,.$
\end{enumerate}
\end{lem}
{\it Proof.} 
In the case that $1\le i\le n-1\,,$ we have 
$s_iw_k=w_ks_i$ for $1\le k\le i-1$ or $2n-i+2\le k\le 2n\,,$
and $s_iw_k=w_ks_{i+1}$ for $i+2\le k\le 2n-i-1\,.$ 
These lead to the desired equalities in (a).\par\smallskip
In the same way, the relations 
$s_nw_k=w_ks_n \quad (k\ne n, n+1)$ and 
$s_\theta w_k=w_k(s_2\cdots s_{n-1})(s_n\cdots s_2) \quad (k\ne 1,2n)$
lead to the relations in (b) and (c). \qquad\qed\par\medskip

Next we consider the action of $R_{\alpha_i}$  on the $\overline{h}_{w_k}\,:$ 
\begin{lem}\label{lem4.5}
\begin{enumerate}
\renewcommand{\labelenumi}{(\alph{enumi})}
\item If $1\le i\le n-1,$ 
$R_{\alpha_i}\overline{h}_{w_k}=\overline{h}_{w_k}$
for each $1\le k\le 2n$ such that $k\ne i,\,i+1,\,2n-i,\,2n-i+1.$
\item$R_{\alpha_n}\overline{h}_{w_k}=\overline{h}_{w_k}$
for each $1\le k\le 2n$ such that
$k\ne n,\,n+1.$
\item $R_{\delta-\theta}\overline{h}_{w_k}=\overline{h}_{w_k}$
for each $2\le k\le 2n-1\,.$ 
\end{enumerate}
\end{lem}
{\it Proof.} Since $w_k^{-1}\alpha_i=\alpha_i>0$  for  $ 1\le k\le i-1$
(then $i\ge 2$)\,,
we have
\begin{align*}
R_{\alpha_i}h_{w_k}&=a_{\alpha_i}h_{w_k}+b_{\alpha_i}h_{s_iw_k}\,,\\
R_{\alpha_i}h_{s_iw_k}&=c_{\alpha_i}h_{s_iw_k}+d_{\alpha_i}h_{w_k}\,.
\end{align*}
These imply 
$$R_{\alpha_i}(h_{w_k}+h_{s_iw_k})=h_{w_k}+h_{s_iw_k}\,,$$
following from the relations
$a_{\alpha_i}+d_{\alpha_i}=b_{\alpha_i}+c_{\alpha_i}=1\,.$
Hence, noting $s_iw_k=w_ks_{i},$ we obtain 
$R_{\alpha_i}\overline{h}_{w_k}=\overline{h}_{w_k}\,.$
Other cases are similarly derived. 
\quad$\square$\par\medskip
\begin{lem}\label{lem4.6}
\begin{enumerate}
\renewcommand{\labelenumi}{(\alph{enumi})}
\item For $1\le i\le n-1\,,$
\begin{equation*}
\left\{
\begin{aligned}
R_{\alpha_i}\overline{h}_{w_i}&=
a_{\alpha_i}\overline{h}_{w_i}+b_{\alpha_i}\overline{h}_{w_{i+1}}\,,\\
R_{\alpha_i}\overline{h}_{w_{i+1}}&=
c_{\alpha_i}\overline{h}_{w_{i+1}}+d_{\alpha_i}\overline{h}_{w_i}\,,
\end{aligned}
\right.
\,
\left\{
\begin{aligned}
R_{\alpha_i}\overline{h}_{w_{2n-i}}&=
a_{\alpha_i}\overline{h}_{w_{2n-i}}+b_{\alpha_i}\overline{h}_{w_{2n-i+1}}\,,\\
R_{\alpha_i}\overline{h}_{w_{2n-i+1}}&=
c_{\alpha_i}\overline{h}_{w_{2n-i+1}}+d_{\alpha_i}\overline{h}_{w_{2n-i}}\,.
\end{aligned}
\right.
\end{equation*}
\item
\begin{equation*}
\left\{\quad
\begin{aligned}
R_{\alpha_n}\overline{h}_{w_{n}}&=
a_{\alpha_n}\overline{h}_{w_{n}}+b_{\alpha_n}\overline{h}_{w_{n+1}}\,,\\
R_{\alpha_n}\overline{h}_{w_{n+1}}&=
c_{\alpha_n}\overline{h}_{w_{n+1}}+d_{\alpha_n}\overline{h}_{w_{n}}\,.
\end{aligned}
\right.
\end{equation*}
\item
\begin{equation*}
\left\{\quad
\begin{aligned}
R_{\delta-\theta}\overline{h}_{w_{2n}}&=
a_{\delta-\theta}\overline{h}_{w_{2n}}+
q^{-\lambda}b_{\delta-\theta}\overline{h}_{w_{1}}\,,\\
R_{\delta-\theta}\overline{h}_{w_{1}}&=
c_{\delta-\theta}\overline{h}_{w_{1}}+
q^{\lambda}d_{\delta-\theta}\overline{h}_{w_{2n}}\,.
\end{aligned}
\right.
\end{equation*}
\end{enumerate}
\end{lem}
{\it Proof.}
This follows almost immediately from the 
definitions.\qquad\qed\par\medskip
At this stage, by combination of the above lemmas, we 
obtain the following:

In case of $1\le i\le n-1\,,$ we have
{\allowdisplaybreaks
\begin{align*}
&r_{s_i}\,\Psi =
\sum_{1\le k\le 2n}\langle\,s_i\varphi_{w_k}\rangle
\overline{h}_{s_iw_k} =
\{\sum
\begin{Sb}
k\ne i,\,i+1, \\[2pt] 2n-i,\, 2n-i+1
\end{Sb}
+\sum
\begin{Sb}
k= i,\,i+1,\\[2pt] 2n-i,\, 2n-i+1
\end{Sb}
\}
\langle\,s_i\varphi_{w_k}\rangle
\overline{h}_{s_iw_k}\\[3pt]
&=\sum
\begin{Sb}
k\ne i,\,i+1, \\[2pt] 2n-i,\, 2n-i+1
\end{Sb}
\langle\,\varphi_{w_k}\rangle
\overline{h}_{w_k}\\[2pt]
&+\left\{b_{\alpha_i}\langle\,\varphi_{w_i}\rangle
+c_{\alpha_i}\langle\,\varphi_{w_{i+1}}\rangle\right\}
\overline{h}_{s_iw_i}
+\left\{a_{\alpha_i}\langle\,\varphi_{w_i}\rangle
+d_{\alpha_i}\langle\,\varphi_{w_{i+1}}\rangle\right\}
\overline{h}_{s_iw_{i+1}}\\[2pt]
&+\left\{b_{\alpha_i}\langle\,\varphi_{w_{2n-i}}\rangle
+c_{\alpha_i}\langle\,\varphi_{w_{2n-i+1}}\rangle\right\}
\overline{h}_{s_iw_{2n-i}}
+\left\{a_{\alpha_i}\langle\,\varphi_{w_{2n-i}}\rangle
+d_{\alpha_i}\langle\,\varphi_{w_{2n-i+1}}\rangle\right\}
\overline{h}_{s_iw_{2n-i+1}}\\[3pt]
&=\sum
\begin{Sb}
k\ne i,\,i+1, \\[2pt] 2n-i,\, 2n-i+1
\end{Sb}
\langle\,\varphi_{w_k}\rangle
\overline{h}_{w_k}\\[2pt]
&+\langle\,\varphi_{w_i}\rangle
\left\{b_{\alpha_i}\overline{h}_{w_{i+1}}
+a_{\alpha_i}\overline{h}_{w_{i}}\right\}
+\langle\,\varphi_{w_{i+1}}\rangle
\left\{c_{\alpha_i}\overline{h}_{w_{i+1}}
+d_{\alpha_i}\overline{h}_{w_{i}}\right\}\\[2pt]
&+\langle\,\varphi_{w_{2n-i}}\rangle
\left\{b_{\alpha_i}\overline{h}_{w_{2n-i+1}}
+a_{\alpha_i}
\overline{h}_{w_{2n-i}}\right\}
+\langle\,\varphi_{w_{2n-i+1}}\rangle
\left\{c_{\alpha_i}
\overline{h}_{w_{2n-i+1}}
+d_{\alpha_i}
\overline{h}_{w_{2n-i+1}}\right\}\\[3pt]
&=R_{\alpha_i}\,\Psi.
\end{align*}}
Similarly, in the case $i=n$, we have
{\allowdisplaybreaks
\begin{align*}
&r_{s_n}\,\Psi =
\biggl\{\;\sum
\begin{Sb}
1\le k\le 2n\\[2pt]
k\ne n,\, n+1
\end{Sb}
+\sum_{k=n,\, n+1}\;\biggr\}
\langle\,s_n\varphi_{w_k}\rangle
\overline{h}_{s_nw_k}\\[3pt]
&=\sum_{k\ne n,n+1}\langle\,\varphi_{w_k}\rangle
\overline{h}_{w_k}\\[2pt]
&+\left\{b_{\alpha_n}\langle\,\varphi_{w_n}\rangle
+c_{\alpha_n}\langle\,\varphi_{w_{n+1}}\rangle\right\}
\overline{h}_{s_nw_n}
+\left\{a_{\alpha_n}\langle\,\varphi_{w_n}\rangle
+d_{\alpha_n}\langle\,\varphi_{w_{n+1}}\rangle\right\}
\overline{h}_{s_nw_{n+1}}\\[3pt]
&=\sum_{k\ne n,n+1}\langle\,\varphi_{w_k}\rangle
\overline{h}_{w_k}\\[2pt]
&+\langle\,\varphi_{w_n}\rangle\left\{b_{\alpha_n}
\overline{h}_{w_{n+1}}
+a_{\alpha_n}\overline{h}_{w_{n}}\right\}
+\langle\,\varphi_{w_{n+1}}\rangle
\left\{c_{\alpha_n}\overline{h}_{w_{n+1}}
+d_{\alpha_n}\overline{h}_{w_{n}}\right\}\\[3pt]
&=R_{\alpha_n}\,\Psi.
\end{align*}}
Finally, if $i=0\,,$  by noting that
$$s_0\overline{h}_{w_{1}}=q^\lambda\overline{h}_{s_\theta w_{1}}
\quad\text{and}\quad
s_0\overline{h}_{w_{2n}}=q^{-\lambda}\overline{h}_{s_\theta w_{2n}}\,,$$
we have
{\allowdisplaybreaks
\begin{align*}
r_{s_0}\,\Psi 
&=\sum
\begin{Sb}
1\le k\le 2n\\[2pt]
k\ne 1,\,2n
\end{Sb}
\langle\,\varphi_{w_k}\,\rangle
\overline{h}_{s_\theta w_k}
+\langle\,s_0\varphi_{w_1}\rangle q^\lambda\overline{h}_{s_\theta w_1}
+\langle\,s_0\varphi_{w_{2n}}\rangle 
q^{-\lambda}\overline{h}_{s_\theta w_{2n}}\\[3pt]
&=\sum_{k\ne 1,\,2n}\langle\,\varphi_{w_k}\,\rangle
\overline{h}_{w_k}
+\langle\,s_0\varphi_{w_1}\,\rangle q^\lambda\overline{h}_{w_{2n}}
+\langle\,s_0\varphi_{w_{2n}}\,\rangle q^{-\lambda}\overline{h}_{w_{1}}\\[3pt]
&=\sum_{k\ne 1,\,2n}\langle\,\varphi_{w_k}\,\rangle\overline{h}_{w_k}
+\left\{
c_{\delta-\theta}\langle\,\varphi_{w_{1}}\,\rangle
+q^{-\lambda}b_{\delta-\theta}\langle\,\varphi_{w_{2n}}\,\rangle
\right\}\overline{h}_{w_1}\\[2pt]
&\quad+\left\{
q^{\lambda}d_{\delta-\theta}\langle\,\varphi_{w_{1}}\,\rangle
+a_{\delta-\theta}\langle\,\varphi_{w_{2n}}\,\rangle
\right\}
\overline{h}_{w_{2n}}\\[3pt]
&=\sum_{k\ne 1,\,2n}\langle\,\varphi_{w_k}\,\rangle
\overline{h}_{w_k}
+\langle\,\varphi_{w_1}\,\rangle 
\left\{
c_{\delta-\theta}\overline{h}_{w_{1}}
+q^{\lambda} d_{\delta-\theta}
\overline{h}_{w_{2n}}\,\right\}\\[2pt]
&\quad+\langle\,\varphi_{w_{2n}}\,\rangle
\left\{
a_{\delta-\theta}\overline{h}_{w_{2n}}
+q^{-\lambda}b_{\delta-\theta}\overline{h}_{w_{1}}
\right\}\\[3pt]
&=R_{\delta-\theta}\,\Psi.
\end{align*}}
This completes the proof of Proposition 3.1.\qquad\qed\par\medskip
\section{Macdonald polynomials}
Macdonald introduced the $q$-difference operators \cite{Mac1}
to define his orthogonal polynomials
associted with root sytems. 
In the case of a root system of  type $C_n\,,$ the 
$q$-difference operator to define such a polynomial 
is given by
$$E=\sum_{a_1,\ldots, a_n=\pm 1}\;\prod_{1\le i<j\le n}
\frac{1-ty_i^{a_i}y_j^{a_j}}{1-y_i^{a_i}y_j^{a_j}}
\prod_{1\le i\le n}
\frac{1-ty_i^{2a_i}}{1-y_i^{2a_i}}T_{y_i}^{\frac{1}{2}a_i}\,,
$$
where
$$(T_{y_i}f)(y_1,\ldots,y_n)=f(y_1,\ldots,qy_i,\ldots,y_n)\,.$$
Its eigenvalue is known to be
\begin{align*}
c_\mu&=\sum_{a_1,\ldots, a_n=\pm 1}\;
\prod_{j=1}^n q^{\frac{1}{2}\lambda_ja_j}t^{\frac{1}{2}(n-j+1)a_j}\\
&=q^{-\frac{1}{2}(\lambda_1+\cdots+\lambda_n)}
\prod_{j=1}^n(1+t^jq^{\lambda_{n-j+1}})
\end{align*}
with the parameter $\mu=(\lambda_1,\ldots,\lambda_n)$\; (We consider only
the special case corresponding to the condition $t_1=t_2=t$). \par\medskip

As for the eigenfunction of the operator $E,$ we easily find the following:
\begin{cor}\label{cor5.1}
The sum 
\begin{equation}
\sum_{i=1}^{2n}t^{i-1}\langle\,\varphi_{w_i}\,\rangle
\end{equation} 
is a solution of the equation attached to the 
parameter $(\lambda, 0,\ldots,0):$
\begin{equation}
 E\psi=c_{(\lambda, 0,\ldots,0)}\psi\,. 
\end{equation}
\end{cor}
{\it Proof.} This is proven by applying the result of 
Kato (Theorem 4.6 in \cite{Kato}) to our Theorem 3.2.
\qquad\qed\par\medskip

We next proceed to simplify the sum (5.1).\smallskip

We  note the equality
\begin{equation}
t^{2n}\prod_{j=1}^n
\frac{\displaystyle\biggl(1-\;\frac{y_j}{x}\biggr)\biggl(1-\;\frac{y_{j}}{x}\biggr)}
{\displaystyle\biggl(1-t\frac{y_j}{x}\biggr)\biggl(1-t\frac{y_{j}^{-1}}{x}\biggr)}
=1+(t-1)\biggl\{\;\sum_{j=1}^{2n}t^{i-1}\varphi_{w_i}\;\biggr\}\,,
\end{equation}
which is demonstrated by using the partial fractions.\par\smallskip
On the other hand, we have 
\begin{equation}
\langle \;
\prod_{j=1}^n
\frac{\displaystyle\biggl(1-\;\frac{y_j}{x}\biggr)\biggl(1-\;\frac{y_{j}}{x}\biggr)}
{\displaystyle\biggl(1-t\frac{y_j}{x}\biggr)\biggl(1-t\frac{y_{j}^{-1}}{x}\biggr)}
\;\rangle=q^\lambda\;\int_{\cal C}\Phi=q^\lambda\;\langle\,1\,\rangle\,,
\end{equation}
which is demonstrated by changing the integration variable 
such that $x\mapsto qx\,.$
Hence, combination of (5.3) and (5.4) gives the relation
$$\sum_{j=1}^{2n}t^{i-1}\langle \;\varphi_{w_i}\;\rangle
=\frac{1-q^{\lambda}t^{2n}}{1-t}\int_{\cal C}\Phi\,.$$

Therefore we reach
\begin{prop}
The function $\int_{\cal C}\Phi$ is a solution to the equation (5.2).
\end{prop}
It should be remarked that this is valid for arbitrary cycle ${\cal C}$
and that linearly independent solutions are obtained by choosing several
cycles. This situation is similar to that studied in \cite{Mi2}.
\par\medskip

In case that the parameter $\mu$ is from the set of partitions, 
the eigenfunction of the form
$$P_\mu (y|q,t)=m_\mu+\sum_{\nu<\mu}a_{\mu\;\nu }m_\nu\,,$$
is the Macdonald polynomial for the root system $C_n.$
Here $m_\mu=\sum_{\nu\in W\mu}e^\nu\,,$  and
$\nu<\mu$ is defined to be $\mu-\nu\in Q^+$ with $Q^+$  
the positive cone of the root lattice.\par\medskip

In our case, to get the Macdonald polynomial, 
it is enough to consider the
case that $\lambda$ is a positive integer and
take the cycle, with the counterclockwise direction,
which encircles the sequence of poles such that
$y_i, y_iq, y_iq^2, \ldots,$ for $1\le i\le n$ and
$y_i^{-1}, y_i^{-1}q, y_i^{-1}q^2, \ldots,$ for $1\le i\le n.$
This is an integral representaion of the Macdonald polynomial
$P_{(\lambda,0,\ldots,0)}(y|q,t)\,.$\par\medskip
  
Moreover, applying the $q$-binomial theorem 
$$\sum_{m\ge 0}\frac{(a)_m}{(q)_i}z^m
=\frac{(az)_\infty}{(z)_\infty}\quad(|z|<1),
\qquad (a)_m=\prod_{0\le k\le m-1}(1-aq^k)$$
and the residue calculus to our integral, we obtain 
an exact expression of the Macdonald polynomial for 
the root system $C_n.$\par\smallskip

\begin{thm}
$$P_{(\lambda,0,\ldots,0)}(y|q,t)=\frac{(q)_\lambda}{(t)_\lambda}\sum
\begin{Sb}
i_1+\cdots+i_{2n}=\lambda\\[3pt]
i_1,\,\ldots,\,i_{2n}\ge 0
\end{Sb}
\frac{(t)_{i_1}\cdots\,(t)_{i_{2n}}}{(q)_{i_1}\cdots\,(q)_{i_{2n}}}
y_1^{i_1-i_{2n}}y_2^{i_2-i_{2n-1}}\cdots\,y_n^{i_n-i_{n+1}}\,.
$$ 
\end{thm}

{\it Remark.} We also have a dirct way to obtain 
the integral representation of the eigenfunction
for (5.2). This will appear in a future paper. 
For the related work, we also refer the reader
to \cite{MN1}

\section{Final comment}
We finally make a comment on the meaning of our elements
$\varphi_{w_i}$ from the viewpoint of the Hecke
algebra. 

Set 
$$T_i=t+\frac{1-te^{\alpha_i}}{1-e^{\alpha_i}}(s_i-1),
\qquad\text{for}\quad 1\le i\le n\,,$$
where $\alpha_i$ is an element of the simple roots 
and $s_i$ a corresponding
generator of the Weyl group $W\,.$ 
This is the Lusztig operator associated with the root system $C_n$ 
(in the special case $t_1=t_2=t$), which satisfies the following:
\begin{align*}
&(T_i-t)(T_i+1)=0 \qquad(1\le i\le n)\,,\\
&T_iT_{i+1}T_i=T_{i+1}T_iT_{i+1} \qquad(1\le i\le n-2)\,,\\
&T_{n-1}T_{n}T_{n-1}T_{n}=T_{n}T_{n-1}T_{n}T_{n-1}\,,\\
&T_iT_{j}=T_jT_{i}\qquad(|i-j|>2)\,.
\end{align*}
These are the fundamental relations for the Hecke algebra  
$H(W)$ associated with
the root system of type $C_n\,.$
The action of the Lusztig operator on our $\varphi_{w_{i}}$ 
is given as follows.
\begin{prop}For $1\le k\le n\,;$
\begin{align*}
&\begin{cases}
&T_i\varphi_{w_{k}}=t\varphi_{w_{k}}\,,\qquad 
i\ne k-1,\, k\,,\\
&T_{k-1}\varphi_{w_{k}}=(t-1)\varphi_{w_{k}}+\varphi_{w_{k-1}}\,,\\
&T_k\varphi_{w_{k}}=t\varphi_{w_{k+1}}\,,
\end{cases}\\
&\begin{cases}
&T_i\varphi_{w_{n+k}}=t\varphi_{w_{n+k}}\,,\qquad 
i\ne n-k,\, n-k+1\,,\\
&T_{n-k+1}\varphi_{w_{n+k}}
=(t-1)\varphi_{w_{n+k}}+\varphi_{w_{n+k-1}}\,,\\
&T_{n-k}\varphi_{w_{n+k}}=t\varphi_{w_{n+k+1}}\,.
\end{cases}
\end{align*}
\end{prop}
This shows that the vector space 
$\oplus_{i=1}^{2n}{\Bbb C}\varphi_{w_{i}}$
gives the representaion of the Hecke algebra 
$H(W)$ for the $C_n$ type.
Moreover, we can also obtain the representation of the 
affine Hecke algebra in the space of the $q$ de Rham cohomology. 
See \cite{MN1} for $A_{n-1}$ case.\par\smallskip 
In any case, we expect that such a basis attached to
the action of the Hecke algebras could be generalized to the case 
of higher representaions. This is our future problem.\par\medskip

{\it Acknowledgement.} The author wishes to thank Professor Shin-ichi Kato
for valuable suggestion.

\bigskip
\begin{flushleft}
\begin{sc}
Katsuhisa Mimachi\\
Department of Mathematics\\
Kyushu University 33\\
Hakozaki, Fukuoka 812-81\\
Japan\\
\end{sc}
\smallskip 
{\it E-mail address\/}: mimachi{\char'100}
math.kyushu-u.ac.jp\\
\end{flushleft}
\end{document}